\documentclass[11pt]{amsart}
\usepackage{amsmath, amssymb, amsfonts, amsthm, graphicx, eucal, overpic}
\usepackage{color, url}
\newcommand{\fhi}{\varphi}
\newcommand{\numbersystem}[1]{\mathbb{#1}}
\newcommand{\D}{\numbersystem{D}}
\newcommand{\R}{\numbersystem{R}}
\newcommand{\K}{\numbersystem{K}}
\newcommand{\Q}{\numbersystem{Q}}
\newcommand{\C}{\numbersystem{C}}
\newcommand{\Z}{\numbersystem{Z}}
\newcommand{\F}{\numbersystem{F}}
\newcommand{\HQ}{\numbersystem{H}}
\newcommand{\kar}[1]{\chi(#1)}
\newcommand{\collection}[1]{{\mathcal#1}}
\newcommand{\setS}{\collection{S}}
\newcommand{\setP}{\collection{P}}
\newcommand{\setB}{\collection{B}}
\newcommand{\setR}{\collection{R}}
\newcommand{\setL}{\collection{L}}

\newcommand{\CS}{\collection{S}}

\newcommand{\card}[1]{\lvert#1\rvert}

\newcommand{\lift}[1]{\raisebox{1.5ex}[0pt]{#1}}

\theoremstyle{plain}
\newtheorem*{sgtheorem}{SG Theorem}
\newtheorem*{mrtheorem}{MR Theorem}

\newtheorem*{ktheorem}{Kelly's Theorem}
\newtheorem*{qtheorem}{Quaternion Theorem}
\newtheorem*{hlemma}{Hirzebruch Lemma}

\newtheorem{theorem}{Theorem}

\newtheorem{primetheorem}{Theorem}

\newtheorem{lemma}{Lemma}
\newtheorem{corollary}[lemma]{Corollary}
\newtheorem{proposition}[lemma]{Proposition}

\newtheorem{problem}{Problem}

\begin{document}

\title{The Sylvester-Gallai Theorem, colourings and algebra}
\thanks{This material is based upon work supported by the National Research Foundation under Grant number 2053752.}
\author{Lou M. Pretorius}
\address{Department of Mathematics and Applied Mathematics,
        University of Pretoria,
        Pretoria 0002, South Africa}
\email{\texttt{lpretor@scientia.up.ac.za}}
\author{Konrad J. Swanepoel}
\address{Department of Mathematical Sciences,
        University of South Africa,
        PO Box 392, UNISA 0003, South Africa}
\email{\texttt{swanekj@unisa.ac.za}}
%\date{\today}
\subjclass[2000]{Primary 51A45, Secondary 05B25, 51E21.}

\begin{abstract}
Our point of departure is the following simple common generalisation of the Sylvester-Gallai theorem and the Motzkin-Rabin theorem:
\begin{quote}
\em Let $\setS$ be a finite set of points in the plane, with each point coloured red or blue or with both colours.
Suppose that for any two distinct points $A,B\in\setS$ sharing a colour there is a third point $C\in\setS$, of the other colour, collinear with $A$ and $B$.
Then all the points in $\setS$ are collinear.
\end{quote}
We define a chromatic geometry to be a simple matroid for which each point is coloured red or blue or with both colours, such that for any two distinct points $A,B\in\setS$ sharing a colour there is a third point $C\in\setS$, of the other colour, collinear with $A$ and $B$.
This is a common generalisation of proper finite linear spaces and properly two-coloured finite linear spaces, with many known properties of both generalising as well.
One such property is Kelly's complex Sylvester-Gallai theorem.
We also consider embeddings of chromatic geometries in Desarguesian projective spaces.
We prove a lower bound of $51$ for the number of points in a $3$-dimensional chromatic geometry in projective space over the quaternions.
Finally, we suggest an elementary approach to the corollary of an inequality of Hirzebruch used by Kelly in his proof of the complex Sylvester-Gallai theorem.
\end{abstract}

\maketitle

\section{Introduction}\label{introduction}
\subsection{The Sylvester-Gallai Theorem}\label{firstsubsection}
The cubic curve with homogeneous equation $x^3+y^3+z^3+xyz=0$ in the complex 
projective plane has nine inflection points given by homogeneous coordinates
\[ \setS = \bigcup_{\omega^3=1} \{ (0,-1,\omega), (\omega, 0, -1), (-1,\omega, 0) \}. \]
The set $\setS$ is non-collinear, yet the line through any two points of the set contains a third point of the set.
This property of cubic curves goes back to Pl\"ucker's book of 1835 \cite{Plucker}; see also \cite{DeVries}.
In 1893 Sylvester \cite{Sylvester} posed a question where he essentially asked for a proof that there is no such finite non-collinear set in the real plane.
Erd\H os independently discovered this question in 1933 and Gallai was the first to solve it \cite{MR84e:01073}, but his solution appeared in print only in 1944 \cite{Erdos, Steinberg}.
We next state his theorem, usually called the Sylvester-Gallai theorem.

\begin{sgtheorem}
Let $\setS$ be a finite set of points in the plane.
Suppose that for any two distinct points $A,B\in\setS$ there is a third point $C\in\setS$ collinear 
with $A$ and $B$.
Then all the points in $\setS$ are collinear.
\end{sgtheorem}
The plane in the SG Theorem can be chosen to be the real affine plane or the real 
projective plane.
The first published proof (of the dual theorem) is due to Melchior \cite{MR3:13c} and uses Euler's polyhedral formula.
Subsequently many proofs appeared, as well as extensions and generalisations \cite{MR46:6148, MR92b:52010, MR96m:52025}.
The best known proof is probably Kelly's proof using the least non-zero perpendicular distance between a point $P\in\setS$ and a line through two other points of $\setS$.
This proof was first published in \cite{Coxeter} (see also \cite{MR2001j:00001} for this proof ``from the book'').

\subsection{Colourings}
Another proof of the SG Theorem is due to Motzkin \cite{MR12:849c}, its simplicity rivalling that of Kelly's famous proof.
It proves the dual statement.
Motzkin used the same idea to prove the following chromatic version of the SG Theorem, which was originally a question posed by Ronald Graham \cite{MR2000e:05036}.
\begin{mrtheorem}
Let $\setS$ be a finite set of points in the plane, each coloured red or blue.
Suppose that for any two distinct points $A,B\in\setS$ of the same colour there is a third point $C\in\setS$ of the other colour, collinear with $A$ and $B$.
Then all the points in $\setS$ are collinear.
\end{mrtheorem}
Again the plane can be either the real affine plane or the real projective plane.
This theorem was independently proved by Michael Rabin, and he and Motzkin planned to write a joint paper (cited for example in \cite{MR41:3305} and \cite{KR}), but the paper never appeared \cite{MR2000e:05036}.
Motzkin's proof is written up in the two surveys \cite{MR96m:52025} and \cite{MR92b:52010}.
The first published proof of the MR Theorem is by Chakerian \cite{MR41:3305}, who proves the dual theorem using the Euler formula (see also \cite{EML}).
See \cite{PreSwa} for a self-contained proof from the axioms of ordered geometry without dualising.
It is not known whether this theorem still holds in the plane over the complex numbers (see Problem~\ref{p3a} in Section~\ref{problems}).
We now make the simple observation that the SG Theorem and the MR Theorem have a common generalisation.

\begin{theorem}\label{thm1}
Let $\setS$ be a finite set of points in the plane, with each point coloured red or blue or with both colours.
Suppose that for any two distinct points $A,B\in\setS$ sharing a colour there is a third point $C\in\setS$ of the other colour, collinear with $A$ and $B$.
Then all the points in $\setS$ are collinear.
\end{theorem}
This is clearly a generalisation of the MR Theorem.
The SG Theorem is obtained by colouring each point with both colours.
The following is a reformulation.
\begin{corollary}
Let $\setB$ and $\setR$ be two finite sets in the plane such that $\setB\cup\setR$ is non-collinear.
Then there exists a line $\ell$ such that $\card{\ell\cap\setB\cap\setR}\leq 2\leq\card{\ell\cap(\setB\cup\setR)}$, and $\ell\cap(\setB\cup\setR)$ is contained in either $\setB$ or $\setR$.
\end{corollary}
The MR Theorem corresponds to the case where $\setB$ and $\setR$ are disjoint, and the SG Theorem where $\setB=\setR$.
For convenience we present the proof of Theorem~\ref{thm1}, which is essentially Motzkin's proof of the dual MR Theorem.
We first formulate the dual statement.
\begin{primetheorem}\label{prime}
Let $\setS$ be a finite set of lines in the real projective plane, with each line coloured red or blue or with both colours.
Suppose that for any two distinct lines $\ell, m\in\setS$ sharing a colour there is a third line $n\in\setS$ of the other colour, concurrent with $\ell$ and $m$.
Then all the lines in $\setS$ are concurrent.
\end{primetheorem}

\begin{proof}
Assume that the lines of $\setS$ are not concurrent.
Since there are then at least three lines in $\setS$, there are two lines $a$ and $b$ sharing a colour, say blue.
If all the red lines pass through $C:=a\cap b$, then $a$ and any line not passing through $C$ (which is then necessarily blue) contradict the hypothesis.
Thus there is a red line $c$ not concurrent with $a$ and $b$ (Figure~\ref{fig1}).
\begin{figure}[h]
\setlength{\unitlength}{2cm}
\begin{center}
\begin{picture}(2.5,1.6)(-0.3,-0.3)
\color{blue}\put(0.8,1.2){\line(1,-1){1.3}} \normalcolor \put(1.3,0.75){$a$}
\color{blue}\put(1.2,1.2){\line(-1,-1){1.3}} \normalcolor \put(0.3,0.5){$b$}
\color{blue}\put(0.7,-0.3){\line(1,1){1.2}} \normalcolor \put(1.95,0.75){$e$}
\thicklines
\color{red}\put(-0.3,0){\line(1,0){2.7}} \normalcolor \put(2.25,0.05){$c$}
\color{red}\put(1,1.2){\line(0,-1){1.45}} \normalcolor \put(0.85,0.4){$d$}
\thinlines
\put(0,0){\circle*{0.05}} \put(-0.05,-0.2){$A$}
\put(2,0){\circle*{0.05}} \put(1.85,-0.2){$B$}
\put(1,1){\circle*{0.05}} \put(0.75,0.93){$C$}
\put(1,0){\circle*{0.05}} \put(1.05,-0.2){$D$}
\put(1.5,0.5){\circle*{0.05}} \put(1.6,0.45){$E$}
\end{picture}
\end{center}
\caption{Proof of Theorem~\ref{prime}}\label{fig1}
\end{figure}
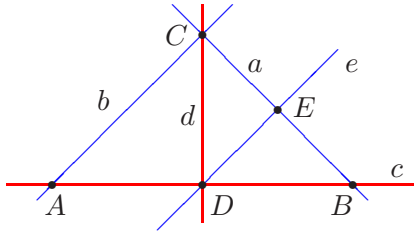
By hypothesis there is a red line $d$ through $C$.
Choose the line at infinity $\omega\neq c$ such that $\omega\cap c$ and $D:=c\cap d$ are separated by $A:=b\cap c$ and $B:=a\cap c$.
Then $D$ is between $A$ and $B$ in the real affine plane obtained by removing $\omega$ and its incident points.
The remainder of the proof plays out in this affine plane.

There is a blue line $e\neq c,d$, passing through $D$, by hypothesis.
By the axiom of Pasch $e$ intersects one of the edges $AC$ or $BC$.
Assume without loss of generality that $e$ intersects $a$ in $E$ between $B$ and $C$.
By renaming $A'=B, B'=C, C'=D, D'=E$, we obtain the same situation as before, only in a smaller (in the sense of containment) triangle $\triangle A'B'C'$ and with the colours switched.
Thus there is a line $e'$ through $D'$, and we may apply the axiom of Pasch again.

Since there are only finitely many lines in $\setS$, this process cannot continue indefinitely.
We thus eventually obtain two lines of the same colour not concurrent with a third of the other colour, contradicting the hypothesis.

It follows that $\setS$ is concurrent.
\end{proof}

\subsection{Algebra}
Since the origin of the Sylvester-Gallai theorem lies in the field of complex numbers, it makes sense to consider the influence of various ground fields.
The first such result, by Kelly \cite{MR87k:14047}, answered a question of Serre \cite{Serre}.
\begin{ktheorem}
Let $\setS$ be a finite set of points in complex projective $n$-space.
Suppose that for any two distinct points $A,B\in\setS$ there is a third point $C\in\setS$ collinear with $A$ and $B$.
Then all the points in $\setS$ are coplanar.
\end{ktheorem}
We generalise Kelly's theorem to a chromatic version where points are allowed to be bicoloured.
\begin{theorem}\label{thm2}
Let $\setS$ be a finite set of points in complex projective $n$-space, with each point coloured red or blue or with both colours.
Suppose that for any two distinct points $A,B\in\setS$ sharing a colour there is a third point $C\in\setS$ of the other colour, collinear with $A$ and $B$.
Then all the points in $\setS$ are coplanar.
\end{theorem}
Theorems~\ref{thm1} and \ref{thm2} suggest that we consider abstract geometries with each point coloured, red, blue, or both, such that for any two points sharing a colour, there is a third point on the line through them of the other colour.
We call such geometries \emph{chromatic geometries} (see Section~\ref{notation} for the exact definition).

In the proof of his theorem, Kelly used the following corollary of a deep inequality of Hirzebruch \cite{MR84m:14037}.
\begin{hlemma}
Let $\setS$ be a non-collinear finite set of points in the complex projective plane.
Then there exists a line $\ell$ such that $2\leq\card{\ell\cap\setS}\leq 3$.
\end{hlemma}
The proof of Theorem~\ref{thm2} also uses this lemma.
No elementary proof of the Hirzebruch Lemma is known.
In Section~\ref{finite} we suggest an elementary approach using the compactness theorem of first order logic and considerations from finite geometry.
However, Kelly's theorem has an elementary proof, which was extended in \cite{EPS} to yield the following.
\begin{qtheorem}
Let $\setS$ be a finite set of points in projective $n$-space over the quaternions.
Suppose that for any two distinct points $A,B\in\setS$ there is a third point $C\in\setS$ collinear 
with $A$ and $B$.
Then all the points in $\setS$ are contained in a three-dimensional flat.
\end{qtheorem}
It is not known whether such a three-dimensional set of points exists.
We show the following lower bounds for the more general, chromatic case, and for arbitrary division rings.
\begin{theorem}\label{thm3}
Let $\setS$ be a finite set of points in projective $n$-space over a division ring $\D$, with each point coloured red or blue or with both colours.
Suppose that for any two distinct points $A,B\in\setS$ sharing a colour there is a third point $C\in\setS$, of the other colour, collinear with $A$ and $B$.
If $\setS$ is not coplanar, then
\begin{itemize}
\item $\card{\setS}\geq 15$, with equality if and only if the characteristic of $\D$ is $2$ and $\setS$ is isomorphic to the $3$-dimensional projective space of order $2$,
\item $\card{\setS}\geq 27$ if the characteristic of $\D$ is not $2$, with equality if and only if the characteristic of $\D$ is $3$ and $\setS$ is isomorphic to the $3$-dimensional affine space of order $3$,
\item $\card{\setS}\geq 51$ if the characteristic of $\D$ is not $2$ or $3$.
\end{itemize}
\end{theorem}
We do not believe that the lower bound of $51$ is sharp.
Note, however, that the set of all the points on three parallel planes in the affine $3$-space over the field of $5$ elements gives an example with $75$ points.

In the next section we introduce the terminology to be used in the remainder of the paper.
In Section~\ref{basic} we consider basic properties of chromatic geometries.
Sections~\ref{proofthm2} and \ref{proofthm3} contain the proofs of Theorems~\ref{thm2} and \ref{thm3}.
In Section~\ref{finite} we consider finite forms of the Hirzebruch Lemma, Kelly's theorem and Theorem~\ref{thm2}.
In Section~\ref{problems} we pose some natural open problems.

\section{Notation and terminology}\label{notation}
\subsection{Geometries and colourings}
We define a \emph{\textup{(}finite incidence\textup{)} geometry} $\setS$ to be a non-empty finite set $\setP$ of \emph{points} and a collection $\setL$ of \emph{lines}, with each line a subset of $\setP$, such that any two points determine a unique line, and any line passes through at least two points.
Three equivalent terms in the literature are \emph{pairwise balanced design} \cite[Part III]{Handbook}, \emph{finite linear space} \cite{BB}, and \emph{finite simple matroid or geometric lattice of rank at most $3$} \cite{Aigner}.
In fact we may consider a matroid of any rank by ignoring flats of dimension higher than $1$.

The \emph{degree} of a point $P$ in a geometry $\setS$ is the number of lines of $\setS$ incident with $P$, and the \emph{size} of a line $\ell$ of $\setS$ is the number of points of $\setS$ incident with $\ell$.

An \emph{SG geometry} is a geometry with at least $3$ points on each line.
Other terms used are \emph{proper finite linear space} \cite{Handbook} and \emph{Sylvester-Gallai design} \cite{MR47:3207}.
An \emph{MR geometry} is a geometry with each point coloured red or blue such that each line contains a red point and a blue point \cite{PreSwa2}.
We call such a colouring a (proper) \emph{$2$-colouring} of the geometry.
An MR geometry is equivalent to a finite linear space which has a blocking set, with the blocking set singled out.
A \emph{chromatic geometry} is a geometry with each point coloured red or blue, or with both colours, such that for any two points $A, B$ sharing a colour, the line $AB$ contains a third point $C$ of the other colour.
A point that is both red and blue is \emph{bicoloured}.
It is immediate that an MR geometry is a chromatic geometry.
If we bicolour each point of an SG geometry we also obtain a chromatic geometry.

Let $\setS$ be any finite simple matroid of rank $m$.
The \emph{contraction} (or \emph{residue}) of $\setS$ at a $k$-flat $\setS'$, denoted by $\setS/\setS'$, is the simple matroid of rank $m-k-1$ of which the $d$-flats correspond to the $(d+k+1)$-flats of $\setS$ that contain $\setS'$.
We say that a point $P$ of $\setS$ is \emph{rich} if each line of the contraction $\setS/ P$ has size at least $4$, and \emph{poor} otherwise.

\subsection{Algebra}
By a \emph{division ring} we mean an associative (not necessarily commutative) ring with identity and with each non-zero element invertible.
(If the division ring is commutative, it is of course a \emph{field}).
We denote a division ring by $\D$, a field by $\F$ and the finite field of size $q$ by $\F_q$.
The set of non-zero elements of $\D$ is denoted by $\D^\ast$, and the characteristic of $\D$ by $\kar{\D}$.
As usual, we denote the rational numbers, the algebraic numbers, the real numbers, the complex numbers and the quaternions by $\Q$, $\overline{\Q}$, $\R$, $\C$ and $\HQ$, respectively.
The $n$-dimensional affine and projective spaces over $\D$ are denoted by $\D^n$ and $P^n(\D)$, respectively.
We use homogeneous coordinates for points of $P^n(\D)$.
We also denote $\F_q^n$ by $AG(n,q)$ and $P^n(\F_q)$ by $PG(n,q)$.
For general background on division rings see \cite{Lam} and on projective spaces over a division ring see e.g.\ \cite{Stevenson}.

An \emph{embedding} of a geometry $\setS$ into $P^n(\D)$ is an injection $\fhi: \setS\to P^n(\D)$ such that collinear points are mapped to collinear points, and non-collinear points to non-collinear points.
An embedding of a geometry into $\D^n$ is defined similarly.
We call the embedded geometry $\fhi(\setS)$ a \emph{configuration}.
An embedded SG (MR, chromatic) geometry is called an \emph{SG \textup{(}MR, chromatic\textup{)} configuration}, and abbreviated by SGC, MRC, CC, respectively.
Two embeddings of a geometry into $P^n(\D)$ are \emph{projectively equivalent} if there exists a projective transformation of $P^n(\D)$ mapping the one embedding to the other.
We use the following consequence of Hilbert's Nullstellensatz.
\begin{lemma}\label{HN}
If a geometry $\setS$ is embeddable as a $k$-dimensional subset of $P^n(\F)$, where $\F$ is a field, then $\setS$ is embeddable as a $k$-dimensional subset of $P^n(\F')$ for some finite extension $\F'$ of the prime subfield of $\F$.
\end{lemma}
\begin{proof}
The embeddability of $\setS$, as well as its $k$-dimensionality, can be expressed as a finite set of polynomial equations and inequations, where the coefficients of the polynomials are in the prime subfield $\K$ of $\F$.
Hilbert's Nullstellensatz now asserts that if this set of equations and inequations has a solution, then it has a solution in the algebraic closure of $\K$.
Since there are only finitely many variables, the solutions in fact lie in a finite algebraic extension of $\K$.
\end{proof}

The \emph{dimension} of a subset $\setS$ of $P^n(\D)$ is the dimension of the flat generated by $\setS$.
An $m$-dimensional subset $\setS$ can be considered to be a matroid of rank $m+1$, by means of the structure inherited from $P^n(\D)$, i.e.,
a \emph{$k$-flat} of $\setS$ is the intersection $\setS'$ of $\setS$ with some $k$-flat $\Pi$ of $P^n(\D)$ such that $\setS'$ spans $\Pi$.
If $\setS$ is embedded in $P^n(\D)$ and $\setS'$ is a $k$-flat of $\setS$, then there is a natural associated embedding of the contraction $\setS/\setS'$ in $P^{n-k-1}(\D)$, which can be realised as an image of $\setS$ under a projection of $P^n(\D)$ onto an appropriate $(n-k-1)$-flat of $P^n(\D)$; see \cite[Chapter~VI]{Aigner}.

In Section~\ref{finite} we use some elementary model theory, all of which may be found for example in \cite{Barwise}.

\section{Basic properties of chromatic geometries}\label{basic}
\subsection{Elementary properties}
It is obvious that a flat of a chromatic simple matroid is again a chromatic geometry.
A contraction $\setS/\setS'$ of a chromatic simple matroid $\setS$ at the $k$-flat $\setS'$ is also a chromatic geometry, if we colour each point $P$ of $\setS/\setS'$ with all the colours of the points of $\Pi$ not contained in $\setS'$, where $\Pi$ is the $(k+1)$-flat associated to $P$.
Note that if we take a contraction of an MR geometry, we may end up with a chromatic geometry where some points are bicoloured.
This shows that, when considering contractions of matroids, the notion of a chromatic geometry is more natural than that of an MR geometry.

The following is an easy observation.
\begin{lemma}\label{lemma1}
Let $\setS$ be a chromatic geometry.
Any line through a bicoloured point in $\setS$ has size at least $3$.
Therefore, if a line has only two points, one of the points must be red-only, and the other blue-only.
\end{lemma}

It is a well-known folklore result that Steiner triple systems cannot be two-coloured.
The following is a simple generalisation.
\begin{proposition}
If each line of a non-collinear chromatic geometry has size at most $3$, then it must be a Steiner triple system with each point bicoloured.
\end{proposition}
\begin{proof}
Let $b$, $r$, $p$, denote the number of blue-only, red-only, and bicoloured points, respectively.
Since a line through a blue-only and a bicoloured point must have a further red-only point, and the line through a red-only and a bicoloured point a further blue-only point, it follows that $bp$ and $rp$ count the same lines, hence $bp=rp$.
Let $t$ denote the number of lines of size $2$.

If $p=0$ then we have an MR geometry with all lines of size $2$ or $3$.
It is well-known that this is not possible; the following is a short argument.
All lines of size $3$ have two points of one colour, and one of the other.
It follows that $\binom{b}{2}$ is the total number of blue-blue-red lines, and $\binom{r}{2}$ the total number of red-red-blue lines.
Therefore,
\[ \binom{b+r}{2} = 3\binom{b}{2}+3\binom{r}{2}+t,\]
which simplifies to
\[ (b-r)^2+(b-1)(r-1)=1-t\leq 1,\]
implying that $b,r\leq 2$, easily leading to a contradiction.

Therefore, $p\geq 1$ and $b=r$.
Since the number of lines with three bicoloured points is $\binom{p}{2}/3$, and a line of size $3$ can only be blue-blue-red, red-red-blue, bicoloured-red-blue, or bicoloured, we have as before
\[\binom{b+r+p}{2} = 3\binom{b}{2}+3\binom{r}{2} + 3bp + \binom{p}{2}+t.\]
Using $b=r$ this simplifies to $b(b+p-2)+t= 0$.
If $b\neq 0$ then $b+p-2\leq 0$, giving $b=r=p=1$ and $t=0$, which implies collinearity, a contradiction.
Thus $b=r=t=0$, showing that each line has size $3$ and each point is bicoloured.
\end{proof}

The following is easy to prove.

\begin{lemma}\label{twolines}
A non-collinear chromatic geometry cannot be covered by two lines.
A non-planar chromatic geometry cannot be covered by two $2$-flats.
\end{lemma}

We now characterise the non-collinear chromatic geometries that can be covered by three lines.
This generalises \cite[Lemma~5]{PreSwa2}, asserting that any point of an MR geometry has degree at least $4$.
It is needed in the proofs of Theorems~\ref{thm2} and \ref{thm3}.
\begin{proposition}\label{threelines}
Let $\setS$ be a chromatic geometry covered by three lines.
Then the three lines are either concurrent or pairwise disjoint, and $\setS$ is an SG geometry with each point bicoloured except perhaps the common point of the three lines.
\end{proposition}
\begin{proof}
Let the three lines be $\ell_1, \ell_2, \ell_3$.
Suppose first that the three lines are non-concurrent and not pairwise disjoint.
Then some two of the lines intersect, say $P\in\ell_1\cap\ell_2$, but $P\notin\ell_3$.
For any $Q\in\ell_3\setminus(\ell_1\cup\ell_2)$ the line $PQ$ cannot contain any point other than $P$ or $Q$.
By Lemma~\ref{lemma1} $P$ is not bicoloured, say $P$ is blue-only, and each point of $\ell_3\setminus(\ell_1\cup\ell_2)$ is red-only.
Then $\ell_3$ has to intersect $\ell_1$ or $\ell_2$ in a blue point, say $B=\ell_1\cap\ell_3$.
As above, $B$ is not bicoloured and all points on $\ell_2\setminus(\ell_1\cup\ell_3)$ are red-only.
Choose red points $R_1\in\ell_1$ and $R_2\in\ell_2$.
Then the only point on $R_1R_2$ other than $R_1$ and $R_2$ must be a red-only point on $\ell_3$, a contradiction.

It follows that the three lines are either concurrent or pairwise disjoint.
If they are pairwise disjoint, we may add a new common point to all three lines to make them concurrent, and colour it arbitrarily.
Therefore, we may assume without loss of generality that the three lines have a common point $P$.

Before showing that each point other than $P$ is bicoloured, we first need that each line $\ell_i$ has a blue point other than $P$.
If not, then for some line, say $\ell_3$, we have that all the points in $\ell_3\setminus\{P\}$ are red-only.
Then $P$ is blue-only, and it follows that $\ell_1$ and $\ell_2$ each has a red point, say $R_1\in\ell_1$ and $R_2\in\ell_2$.
Then $R_1R_2$ has to intersect $\ell_3\setminus\{P\}$ in a blue point, a contradiction.
Similarly, each $\ell_i$ has a red point other than $P$.

Suppose now that some point $\neq P$ is not bicoloured, say $B_1\in\ell_1$ is blue-only.
Let $B_2\in\ell_2\setminus\{P\}$ and $B_3\in\ell_3\setminus\{P\}$ be blue.
Let $B_1B_2$ intersect $\ell_3\setminus\{P\}$ in the red point $R_3$ (possibly $R_3=B_3$).
Let $r$ be the number of red points on $\ell_1\setminus\{P\}$, and $b$ the number of blue points on $\ell_2\setminus\{P\}$.
For any red point $R\in\ell_1\setminus\{P\}$, $R_3R$ intersects $\ell_2$ in a blue point $B_R\neq B_2, P$.
Since $R\mapsto B_R$ is an injection, it follows that $r\leq b-1$.
On the other hand, for any blue point $B\in\ell_2\setminus\{P\}$, $B_3B$ intersects $\ell_1$ in a red point $\neq P$.
This gives $b\leq r$, a contradiction.

Finally we show that $\setS$ is an SG geometry.
The only lines for which it is not yet clear that they contain at least $3$ points, are the $\ell_i$.
If they have a common point $P$, then we have already seen that each $\ell_i$ contains a red point and a blue point $\neq P$.
If these red and blue points are different, there are already three points on $\ell_i$.
Otherwise $\ell_i$ contains a bicoloured point, hence must have size at least $3$ by Lemma~\ref{lemma1}.
The only remaining case is if the $\ell_i$ are pairwise disjoint.
Consider e.g.\ $\ell_1$.
Let $B_1$ be a red point and $R_1$ a blue point on $\ell_1$ with $B_1\neq R_1$.
Let $R_2$ be a red point on $\ell_2$.
Then $R_1R_2$ intersects $\ell_3$ in a blue point $B_3$, $B_1B_3$ intersects $\ell_2$ in a red point $S_2\neq R_2$, there is a blue point $B_2\in\ell_2\setminus\{R_2,S_2\}$, and $B_2B_3$ intersects $\ell_1$ in a red point $\neq R_1, B_1$.
This gives $\card{\ell_1}\geq 3$.
\end{proof}
Let $\setS$ be a non-collinear SG geometry covered by three lines $\ell_1, \ell_2, \ell_3$.
If there is a point of concurrency, denote it by $P$.
It is easily seen that each line has the same number of points.
We note in passing that deleting $P$ if necessary, there is a simple correspondence between the isomorphism classes of such SG geometries and the so-called main classes of Latin squares \cite[II.1]{Handbook}.

The following generalises a known property of MR geometries \cite{PreSwa2}.

\begin{proposition}\label{atleast6}
In any chromatic geometry $\setS$
\begin{enumerate}
\item there are at least $6$ blue points and at least $6$ red points, and\label{one}
\item for any line $\ell$ there are at least $3$ blue points and at least $3$ red points not on $\ell$.\label{two}
\end{enumerate}
\end{proposition}
\begin{proof}
We first consider \eqref{two}.
If some point $P\in\ell$ has degree $3$, then by Proposition~\ref{threelines} all points except perhaps $P$ are bicoloured, and each line through $P$ has at least $3$ points.
It follows that there are at least $4$ blue points not on $\ell$, with a similar statement for red points.

Without loss we therefore assume that each point on $\ell$ has degree $\geq 4$.
On each line through a red point on $\ell$ there is at least one blue point.
This gives at least three blue points not on $\ell$.
Similarly, there are at least three red points not on $\ell$.

We now consider \eqref{one}.
If some point has degree $3$, then Proposition~\ref{threelines} provides at least $6$ blue points and at least $6$ red points in $\setS$.
Assume therefore that each point has degree $\geq 4$.
Let $\ell$ be the line through some two blue points.
Then part \eqref{two} gives at least $5$ blue points in $\setS$.
Assume for the sake of contradiction that there are exactly $5$ blue points $B_i$, $i=1,\dots,5$.
Then by part~\eqref{two} no three blue points are collinear, hence all points on $B_iB_j$ other than $B_i$ and $B_j$ are red-only.
For each pair $\{i,j\}$ with $i<j$, choose such a red-only point $R_{ij}\in B_iB_j$.
Then the $R_{ij}$ are all distinct, otherwise one of them would have degree $3$.
Since $R_{12}$ has degree $4$, it must be on a line $\ell$ containing $\geq 3$ other $R_{ij}$ by the pigeon-hole principle.
Then $\card{\ell}\geq 5$.
It follows that all $R_{ij}\in\ell$, since each $R_{ij}$ only has degree $4$.
Let $B_k$ be a blue point on $\ell$.
Since $R_{ik}\in\ell$ for each $i\neq k$, we obtain $B_i\in\ell$ for all $i$, a contradiction.

Therefore, there are at least $6$ blue points, and similarly, at least $6$ red points.
\end{proof}

By a theorem of Bruen \cite{Bruen} the number of points of a single colour in a two-colouring of a finite projective space of order $n$ is at least $n+\sqrt{n}+1$.
In fact it is easy to check that the same proof works for chromatic geometries.
\begin{proposition}
If a projective plane of order $n$ is coloured so that it becomes a chromatic geometry, then the number of red points is at least $n+\sqrt{n}+1$.
\end{proposition}

\subsection{Small chromatic geometries}
The SG geometries on up to $18$ points have been enumerated by Kelly and Nwankpa \cite{MR47:3207}, Brouwer \cite{MR83b:05031}, Heathcote \cite{MR95j:51018}, and Betten and Betten \cite{MR2000f:05020,MR2002g:05042} (Table~\ref{tabl1}).
%In Table~\ref{tabl1} we list the number of non-isomorphic non-collinear SG geometries for each number of points $\leq 18$.
%\begin{table}
%\begin{center}
%\begin{tabular}{|r|r|r|}
%\hline Points & SG & MR \\ \hline
%\hline $\leq 6$ & $0$ & $0$ \\
%\hline $7$ & $1$ & $0$ \\
%\hline $8$ & $0$ & $0$ \\
%\hline $9$ & $1$ & $0$ \\
%\hline $10$ & $1$ & $0$ \\
%\hline $11$ & $1$ & $0$ \\
%\hline $12$ & $3$ & $1$ \\
%\hline $13$ & $7$ & $1$ \\
%\hline $14$ & $1$ & $2$ \\
%\hline $15$ & $119$ & $6$ \\
%\hline $16$ & $398$ & $18$ \\
%\hline $17$ & $161925$ & $82$ \\
%\hline $18$ & $24212890$ & $1000$\\ \hline
%\end{tabular}
%\bigskip
%\end{center}
%\caption{Number of non-collinear SG and MR geometries on up to $18$ points}\label{tabl1}
%\end{table}
\begin{table}
\begin{center}\small
\begin{tabular}{|c|c|c|c|c|c|c|c|c|c|c|c|c|c|}
\hline
Points & $\leq 6$ & 7 & 8 & 9 & 10 & 11 & 12 & 13 & 14 & 15 & 16 & 17 & 18 \\ \hline
SG & 0 & 1 & 0 & 1 & 1 & 1 & 3 & 7 & 1 & 119 & 398 & 161925 & 24212890 \\ \hline
MR & 0 & 0 & 0 & 0 & 0 & 0 & 1 & 1 & 2 & 6   & 18  & 82     & 1000 \\ \hline
\end{tabular}
\smallskip
\end{center}
\caption{Number of non-collinear SG and MR geometries on $\leq 18$ points}\label{tabl1}
\end{table}
Kelly and Nwankpa \cite{MR47:3207} also determined which of the SG geometries on at most $14$ points are embeddable in Desarguesian projective planes or spaces.
The MR geometries on up to $15$ points have been enumerated in \cite{PreSwa2}, and up to $18$ points in \cite{vW} (Table~\ref{tabl1}).
%In Table~\ref{tabl2} we list the number of non-isomorphic non-collinear MR geometries for each number of points $\leq 18$.
%\begin{table}
%\begin{center}
%\begin{tabular}{|c|c|c|c|c|c|c|c|c|}
%\hline
%Points & $\leq 11$ & 12 & 13 & 14 & 15 & 16 & 17 & 18\\
%\hline Number of &&&&&&&&\\
%non-collinear & 0 & 1 & 1 & 2 & 6 & 18 & 82 & 1000\\
%MR geometries &&&&&&&&\\ \hline
%\end{tabular}
%\end{center}
%\caption{}\label{tabl2}
%\end{table}
In the proof of Theorem~\ref{thm3} we need a list of all chromatic geometries up to $10$ points.
\begin{lemma}\label{smallcc}
All non-collinear chromatic geometries on at most $10$ points are listed in Table~\ref{cctable} with necessary and sufficient conditions for embeddability into $P^2(\D)$.
\begin{table}
\begin{center}\small
\begin{tabular}{|c|c|l|l|} \hline
Name & Size & Description & Embeddability \\ \hline\hline
& & Fano plane with each point & \\
\lift{$PG(2,2)$} & \lift{$7$} & bicoloured & \lift{$\kar{\D}=2$} \\ \hline
& & Affine plane of order $3$ & $\kar{\D}=3$ or \\ 
\lift{$AG(2,3)$} &\lift{$9$} & with each point bicoloured & $\exists x\in\D(x\neq 1, x^3=1)$\\ \hline
& & Extension of $AG(2,3)$ by a point & \\
\lift{$AG(2,3)^+$} &\lift{$10$} & at infinity, arbitrarily coloured & \lift{$\kar{\D}=3$} \\ \hline
\end{tabular}
\bigskip
\end{center}
\caption{Non-collinear chromatic geometries on at most $10$ points}\label{cctable}
\end{table}
\end{lemma}
\begin{proof}
Since Table~\ref{cctable} gives all SG geometries up to $10$ points \cite{MR47:3207}, we only have to prove that any chromatic geometry on at most $10$ points is an SG geometry.
Suppose that $\setS$ is a counterexample.
Then $\setS$ has a line of size $2$, and by Proposition~\ref{threelines} each point has degree at least $4$.
If some bicoloured point has degree at least $5$, then there are at least $11$ points by Lemma~\ref{lemma1}.
Therefore, each bicoloured point has degree exactly $4$, and by Lemma~\ref{lemma1} there are at least $9$ points in $\setS$.
By Proposition~\ref{atleast6} there are at least two bicoloured points.
If a line through two bicoloured points has size at least $5$, then there are at least $11$ points by Lemma~\ref{lemma1}.
Therefore, all lines through two bicoloured points have size $3$ or $4$, and we distinguish between two cases.

\boldmath\textbf{The line through some two bicoloured points $A$ and $B$  has size $4$:} \unboldmath
Consider the grid determined by the lines through $A$ and through $B$ (Fig~\ref{fig5}(a)).
\begin{figure}
\begin{center}
\hspace{0.8cm}
\begin{overpic}[scale=0.4]{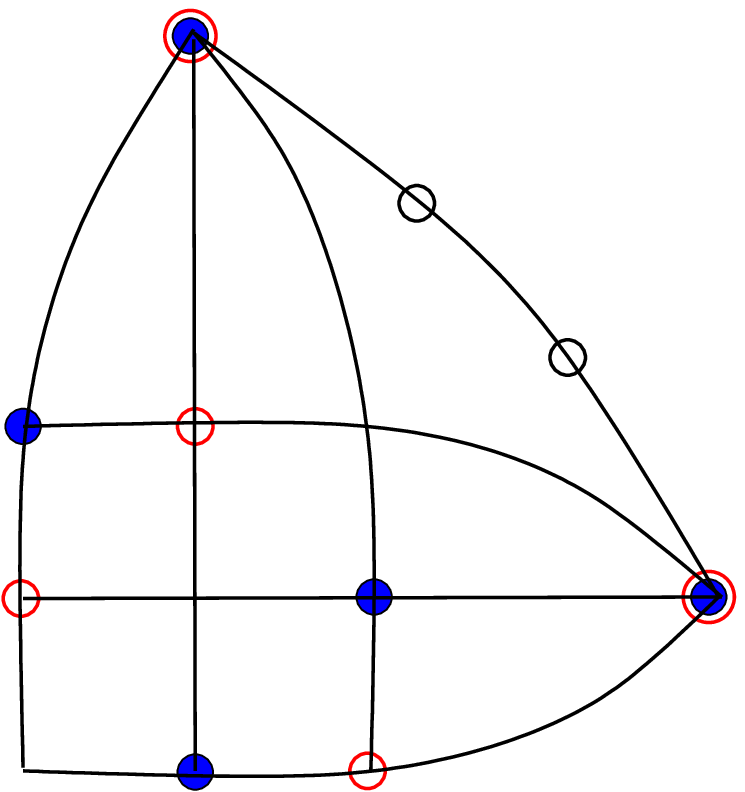}\small
\put(21,-8){$E$}
\put(44,-8){$F$}
\put(89,14){$B$}
\put(50,27){$G$}
\put(26,50){$D$}
\put(27,98){$A$}
\put(-12,45){$C$}
\put(-12,22){$H$}
\put(31,-25){(a)}
\end{overpic}
\hspace{2cm}
\begin{overpic}[scale=0.4]{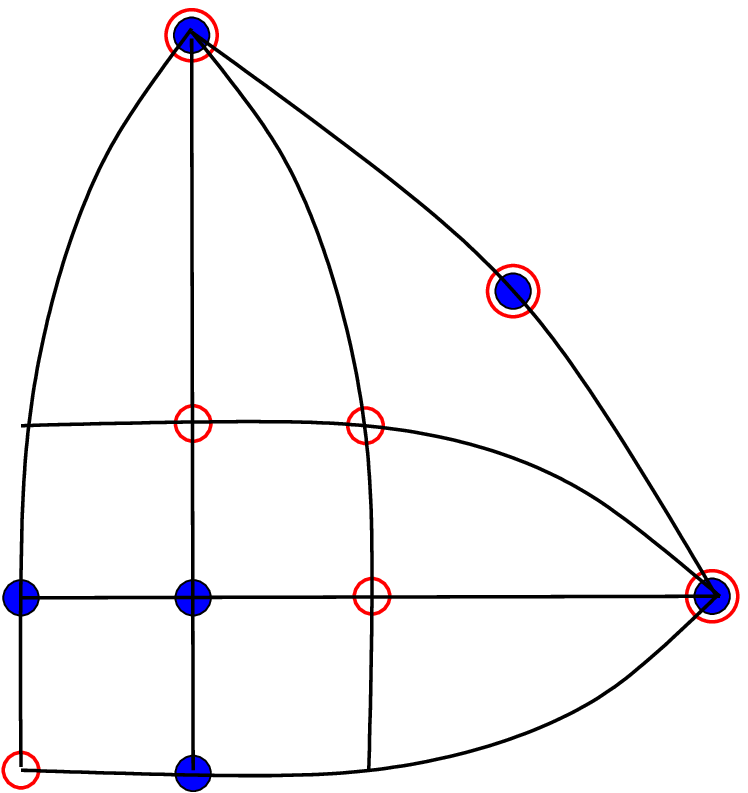}\small
\put(21,-8){$F$}
\put(-10,-6){$D$}
\put(89,14){$C$}
\put(-12,22){$E$}
\put(27,27){$H$}
\put(50,27){$I$}
\put(26,50){$J$}
\put(48,50){$G$}
\put(67,69){$B$}
\put(27,98){$A$}
\put(31,-25){(b)}
\end{overpic}
\bigskip
\end{center}
\caption{Proof of Lemma~\ref{smallcc}}\label{fig5}
\end{figure}
There must be exactly $10$ points, and all lines through $A$ or $B$ except $AB$ must have size $3$.
At least one of the two points on the line of size $2$ must be off $AB$.
Without loss this point $C$, say, is blue-only.
Then $BC$ constains a red-only $D$, $AD$ contains a blue-only $E$, $BE$ contains a red-only $F$, $AF$ contains a blue-only $G$, and $BG$ contains a red-only $H$.
Thus all bicoloured points are on $AB$.
The lines $CE$ and $CG$ contain distinct red points.
These red points have to be on $AB$, and they must be red-only.
Similarly, $DF$ and $DH$ force two blue-only points on $AB$, giving at least $6$ points on $AB$, a contradiction.

\boldmath\textbf{Each line through two bicoloured points has size $3$:} \unboldmath
The bicoloured points form a Steiner triple system.
If this STS is non-collinear, it must be either the Fano plane $PG(2,2)$ or the affine plane $AG(2,3)$.
In the case of $PG(2,2)$ we obtain that there is a unique line through each bicoloured point containing all singly coloured points, a contradiction.
In the case of $AG(2,3)$ we obtain that the four lines through any bicoloured point can contain only bicoloured points, giving that there are no singly coloured points, again a contradiction.

Therefore, the STS is collinear, i.e., there are exactly three bicoloured points $A$, $B$, and $C$, say, on the same line.
Consider the grid determined by the lines through $A$ and the lines through $C$ (Fig.~\ref{fig5}(b)).
Consider a fourth point $D$, say, which is red without loss of generality.
There is a blue $E\in AD$ and a blue $F\in CD$.
The line $EF$ must contain a red-only point $G\neq B$.
Similarly, the line $DG$ contains a blue-only point $H$.
We now obtain red-only points $I\in CH$ and $J\in AH$.
We now have $10$ points, but line $AG$ still needs another blue point, which gives a contradiction.

The embeddability conditions in Table~\ref{cctable} are well-known (see \cite[Theorem~6.1]{Rigby} and \cite[Theorem~3.2 and 3.5]{MR47:3207}) and can all be easily deduced from Propositions~\ref{prop1} and \ref{prop2} below.
\end{proof}

Paul van Wamelen (private communication) has shown that, up to $13$ points, the only chromatic geometry that is not an SG geometry, and that cannot be recoloured to become an MR geometry, is the geometry on $11$ points obtained by bicolouring each point of $AG(2,3)$ and adding two points at infinity (i.e.\ in $PG(2,3)$), one blue-only, and one red-only.

\subsection{Projective embeddings of chromatic geometries covered by three lines}
By Proposition~\ref{threelines}, a chromatic geometry $\setS$ covered by three lines $\ell_1,\ell_2,\ell_3$, is an SG geometry with all points bicoloured, except a possible point of concurrency of the three lines.
We now describe the projective embeddings of such an $\setS$.
Motzkin \cite{MR12:849c} gave an incomplete description.
A correct description was given in the case of fields by Kelly and Nwankpa \cite[Theorem 3.11 and 3.12]{MR47:3207}.
However, this description is still valid for division rings, which is needed in proving Theorem~\ref{thm3}.
\begin{proposition}\label{prop1}
Let $\D$ be a division ring, and let $\setS$ be an SGC in $P^2(\D)$ contained in the union of three concurrent lines intersecting in $P$.
Then $S\setminus\{P\}$ is projectively equivalent to
\[\{(\gamma,0,1), (\gamma,1,1), (-\gamma,1,0) : \gamma\in G\},\]
where $G$ is a finite subgroup of the additive group of $\D$, with $|G|\geq 2$ if $P\in \setS$, and $|G|\geq 3$ if $P\notin \setS$.
\end{proposition}
Note that if we apply this proposition with $G$ consisting of the third roots of unity in $\C$, we do not obtain the set $\CS$ in Section~\ref{firstsubsection}, but rather the set of inflection points of the curve $-x^3+y^3+z^3-xyz=0$, which is clearly projectively equivalent to the original curve.

If $\D$ has a finite additive subgroup $G$, then $\D$ has finite characteristic $p$, and $G$ is isomorphic to the direct sum of finitely many copies of $\Z_p$, the additive group of $\F_p$.
Keeping Proposition~\ref{threelines} in mind, we obtain:

\begin{corollary}\label{deg4}
If $\kar{\D}=0$, then each point in a $2$-dimensional CC has degree at least $4$.

Suppose $\kar{\D}=p$, and $\setS$ is a CC contained in the union of three coplanar concurrent lines $\ell_1, \ell_2, \ell_3$ intersecting in $P$.
Then the number of points in each of the three sets $\ell_i\setminus\{P\}$ is the same positive power of $p$.
\end{corollary}

\begin{proposition}\label{prop2}
Let $\D$ be a division ring, and let $\setS$ be an SGC in $P^2(\D)$ contained in the union of three non-concurrent lines.
Then $S$ is projectively equivalent to
\[\{(1,\gamma,0), (1,0,\gamma), (0,-\gamma,1) : \gamma\in G\},\]
where $G$ is a finite subgroup of the multiplicative group $\D^\ast$ with $|G|\geq 3$.
\end{proposition}
The proofs of Propositions~\ref{prop1} and \ref{prop2} may be found in \cite{PreSwa3}.

\section{Proof of Theorem~\ref{thm2}}\label{proofthm2}
We are now ready to prove Theorem~\ref{thm2}.
Suppose that there exists a CC in $P^n(\C)$ of dimension at least $3$.
Then $P^3(\C)$ contains a CC $\setS$ of dimension $3$.
Consider any point $P\in\setS$, and let $\setS'$ be the contraction of $\setS$ in $P$.
Then $\setS'$ is a CC in $P^2(\C)$.
If some line $\ell$ of $\setS'$ has at most $3$ points, then the corresponding plane through $P$ intersects $\setS$ in a CC where $P$ has degree at most $3$.
This contradicts Corollary~\ref{deg4}.
Therefore, all lines of $\setS'$ have at least $4$ points.
This contradicts the Hirzebruch Lemma, which finishes the proof.
\qed

\section{Proof of Theorem~\ref{thm3}}\label{proofthm3}
We try to mimic the proof of Theorem~\ref{thm2}, where we encountered an SGC with at least $4$ points on each line.
We define a \emph{$k$-SG geometry} to be a finite geometry with at least $k$ points on each line.
A \emph{$k$-SGC} is a $k$-SG geometry embedded into some projective space.
Thus a $3$-SG geometry is an SG geometry, and a $3$-SGC is an SGC.
The contraction of a chromatic configuration in a rich point is a $4$-SGC.
In the absence of an analogue of the Hirzebruch Lemma for division rings we need a list of small $4$-SG geometries.
\begin{lemma}\label{list}
All the non-collinear $4$-SG geometries up to 24 points are given in Table~\ref{4sg} with necessary and sufficient conditions for embeddability into $P^2(\D)$ \textup{(}except for $22.1$--$22.7$, where only necessary conditions are given\textup{)}.
\end{lemma}
\begin{table}
\begin{center}\small
\begin{tabular}{|c|c|l|l|} \hline
Name & Size & Description & Embeddability \\ \hline\hline
$PG(2,3)$ & $13$ & Projective plane of order $3$ & $\kar{\D}=3$ \\ \hline
& & & $\kar{\D}=2$ and \\
\lift{$AG(2,4)$} & \lift{$16$} & \lift{Affine plane of order $4$} & $\exists x\in\D(x\neq 1, x^3=1)$\\ \hline
& & Affine plane of order $4$ & $\kar{\D}=2$ and \\
\lift{$AG(2,4)^+$} & \lift{$17$} & with point at infinity & $\exists x\in\D(x\neq 1, x^3=1)$\\ \hline
& & Punctured projective plane & $\kar{\D}=2$ and \\
\lift{$20.1$} & \lift{$20$} & of order $4$ & $\exists x\in\D(x\neq 1, x^3=1)$\\ \hline
& & Complement of a line & \\
\lift{$20.2$} & \lift{$20$} & of affine plane of order $5$ & \lift{$\kar{\D}=5$} \\ \hline
& & & $\kar{\D}=2$ and \\
\lift{$PG(2,4)$} & \lift{$21$} & \lift{Projective plane of order $4$} & $\exists x\in\D(x\neq 1, x^3=1)$\\ \hline
& & Affine plane of order $5$ with & \\
\lift{$21.1$} & \lift{$21$} & four collinear points removed & \lift{$\kar{\D}=5$} \\ \hline
& & Projective plane of order $5$ & \\
$21.2$ & $21$ & with two lines removed except & $\kar{\D}=5$ \\
&& their point of intersection & \\ \hline
& & Kirkman triple system on $15$ & \\
&& points together with $7$ collinear & $\kar{\D}=2$ is a\\
\lift{$22.1$--$22.7$} & \lift{$22$} &  points, each lying on all lines of a &necessary condition \\
&& parallel class ($7$ geometries) & \\ \hline
& & Projective plane of order $5$ with two & \\
$22.8$ & $22$ & lines removed except for one point & $\kar{\D}=5$ \\
&& (other than the point of intersection) & \\
& & on each & \\ \hline
$24.1$ & $24$ & Punctured affine plane of order $5$ & $\kar{\D}=5$ \\ \hline
& & Projective plane of order $5$ with & \\
&& the points on two lines removed, & \\
\lift{$24.2$} & \lift{$24$} & except for the point of intersection & \lift{$\kar{\D}=5$} \\
&& and three other points on one line & \\ \hline
\end{tabular}
\bigskip
\end{center}
\caption{$4$-SG geometries of $\leq 24$ points}\label{4sg}
\end{table}
\begin{proof}
That this list of $18$ geometries is complete follows from the work of Greig \cite{Greig, Greig2} which corrects and extends Chapter~3 of \cite{BB}.
Since it is not so simple to navigate these general theorems and glean our specific results from them, we give a quick guide through \cite{Greig2}.
Let $\setS$ be a $4$-SG geometry with $\card{\setS}\leq 24$.
We consider different cases, depending on the line sizes that occur in $\setS$.

If some line has size $\geq 7$, then by \cite[Theorem~3.2]{Greig2}, $\setS$ can only be one of $22.1$--$22.7$.

If all lines have the same size $4$, $5$ or $6$, then by \cite[Theorem~6.2]{Greig2}, $\setS$ can only be one of $PG(2,3)$, $PG(2,4)$, or $AG(2,4)$.

If all lines have sizes $4$ or $5$, with both sizes occurring, then by \cite[Theorem~6.9]{Greig2}, $\setS$ can only be one of $20.1$, $AG(2,4)^+$, $24.1$, $20.2$, or $21.1$.

If all lines have sizes $5$ or $6$, with both sizes occurring, then by \cite[Theorem~6.9]{Greig2} there is no possibility for $\setS$.

If all lines have sizes $4$ or $6$, with both sizes occurring, then by \cite[Theorem~6.10]{Greig2}, $\setS$ can only be $21.2$.

If all lines have sizes $4$, $5$, or $6$, with all three sizes occurring, then all points of $\setS$ must have degree $\geq 4$ (since all lines have size $\geq 4$), and $\leq 7$ (otherwise there are $>24$ points).
By \cite[Lemma~6.12]{Greig2}, no point of $\setS$ can have degree $7$.
By \cite[Lemma~6.15]{Greig2}, if at least one point has degree $5$, then $\setS$ can only be one of $22.8$ or $24.2$.
(Note that in the third line of Lemma~6.15 in \cite{Greig2}, the statement ``$v=(n+1)^2$'' should be deleted.)
By \cite[Lemma~6.16]{Greig2}, if all points have degree $6$, then there is no possibility for $\setS$.
At this stage there is still a case that has to be covered, which is seemingly not treated in \cite{Greig2}, namely where all points have degree $4$ or $6$, with degree $4$ occurring.
This case may be finished off as follows.
Let $P$ be a point of degree $4$.
Then all lines not through $P$ must have size $4$.
It then follows easily that all lines through $P$ must have the same size.
However, then only two line sizes occur, a contradiction.

We now consider embeddability.
It is well-known that if $P^2(\D)$ is embedded in $P^2(\D')$, then $\D$ is a subring of $\D'$ \cite[Theorem~8.2.10]{Stevenson}.
This takes care of $PG(2,3)$ and $PG(2,4)$.
It is also known that if the affine plane $\D^2$ is embedded in $P^2(\D')$, and $\D\neq\F_2,\F_3$, then $\D$ is a subring of $\D'$ \cite{Rigby}.
This takes care of $AG(2,4)$ (noting that the condition given is equivalent to $\F_4\subset\D$), as well as its supersets $AG(2,4)^+$, $20.1$, and (again) $PG(2,4)$.

It is shown in \cite{BrCol} that if a set consisting of four parallel lines of a finite affine plane of order $n$ is embeddable in $P^2(\F)$, then $\kar{\F}=n$.
This result carries over to division rings; see \cite{PreSwa3}.
This takes care of $20.2$ and its supersets $21.1$, $21.2$, $22.8$, $24.1$, $24.2$.

By comparing the table of Kirkman triple systems of order $15$ in \cite[I.6.3]{Handbook} with the table of properties of Steiner triple systems of order $15$ in \cite[I.1.2]{Handbook}, we see that each of $22.1$--$22.7$ contains a Fano plane $PG(2,2)$, which forces $\kar{\D}=2$.
\end{proof}

The $4$-SG geometries of size $25$ have been partially classified \cite{Greig}, but we do not know which of them embed into Desarguesian projective planes.
This is the bottleneck in improving Theorem~\ref{thm3}.
Before proving Theorem~\ref{thm3} we note in passing the following two-dimensional analogue.
\begin{corollary}\label{4sgc}
Consider a non-collinear $4$-SGC $\setS$ in $P^n(\D)$.
\begin{itemize}
\item Then $\card{\setS}\geq 13$ with equality if and only if $\setS$ is isomorphic to $PG(2,3)$.
\item If $\kar{\D}\neq 3$ then $\card{\setS}\geq 16$ with equality if and only if $\setS$ is isomorphic to $AG(2,4)$.
\item If $\kar{\D}\neq 2, 3$, then $\card{\setS}\geq 20$ with equality if and only if $\setS$ is isomorphic to the complement of a line in $AG(2,5)$.
\end{itemize}
\end{corollary}

\begin{proof}[Proof of Theorem~\ref{thm3}]

\textbf{Arbitrary characteristic:}
The contraction of $\setS$ at any point is a $2$-dimensional CC, hence has at least $7$ points by Lemma~\ref{smallcc}.
If each line of $\setS$ has size at least $3$, it follows that $\setS$ has at least $(3-1)\times 7+1=15$ points.

Assume then there exists a line $AB$ of $\setS$ of size $2$.
Consider a plane $\Pi$ of $\setS$ passing through $A$ but not through $B$.
Since the degree of $A$ in $\Pi$ is at least $3$, it follows that at least $3$ planes of $\setS$ contain $AB$.
Each of these planes is a CC, hence has at least $7$ points by Lemma~\ref{smallcc}.
It follows that $\setS$ has at least $(7-2)\times 3+2>15$ points.

If $\setS$ has exactly $15$ points, then by the above each line must have size $3$, and each contraction at a point must have exactly $7$ points, hence must be the Fano plane $PG(2,2)$.
It then follows easily that $\setS$ is isomorphic to $PG(3,2)$ \cite[Theorem~7.4.3]{BB}.

\textbf{Characteristic \boldmath$\neq 2$\unboldmath:}
By Lemma~\ref{smallcc}, each plane of $\setS$ has size $\geq 9$.
Suppose some plane $\Pi$ has size $\geq 10$.
Then some point $P\in\Pi$ must have degree $d\geq 4$ in $\Pi$ (otherwise $\card{\Pi}\leq 7$).

Assume that all lines $\ell$ of $\setS$ such that $\ell\cap\Pi=P$ have size $\geq 4$.
If there are $\geq 6$ such lines, we obtain $\card{\setS}\geq 6\times 3+\card{\Pi}\geq 28$.
Thus without loss of generality, there are only $\leq 5$ such lines.
Consider the contraction $\setS/ P$.
The line $\pi$ of $\setS/ P$ corresponding to $\Pi$ has size $d$, and there are $\leq 5$ points of $\setS/ P$ not on $\pi$.
By Proposition~\ref{atleast6}, one of the points not on $\pi$, say $A$, must be bicoloured.
Therefore, all lines in $\setS/ P$ through $A$ have size $\geq 3$ by Lemma~\ref{lemma1}.
Since $A$ has degree $\geq d$ in $\setS/ P$, we obtain $\geq 2d+1$ points in $\setS/ P$, with $\geq d+1$ not on $\pi$.
Therefore, $d+1\leq 5$.
It follows that $d=4$, and $\card{\setS/ P}=9$.
Then $\setS/ P$ must be an $AG(2,3)$ by Lemma~\ref{smallcc}.
However, $AG(2,3)$ does not have a line of size $4$, a contradiction.

Therefore, some line $\ell$ with $\ell\cap\Pi=P$ has size $\leq 3$.
Let the planes through $\ell$ and the $d$ lines through $P$ in $\Pi$ be $\Pi_1,\dots,\Pi_d$.
Each of these planes has size at least $9$.
If $\card{\ell}=2$, we obtain $\card{\setS}\geq(9-2)d+2\geq 7\times 4+2>27$, and if $\card{\ell}=3$, we obtain $\card{\setS}\geq(9-3)d+3\geq 6\times 4+3=27$.
If $\card{\setS}=27$, then $d=4$, and each $\Pi_i$ is an $AG(2,3)$ by Lemma~\ref{smallcc}.
It follows that each line through $P$ in $\Pi$ has size $3$, giving $\card{\Pi}=(3-1)\times 4+1=9$, a contradiction.
Therefore, $\card{\setS}>27$.

We have shown that if some plane of $\setS$ has size $\geq 10$, then $\card{\setS}>27$.
If on the other hand all planes have size $9$, then they are all $AG(2,3)$ by Lemma~\ref{smallcc}.
This implies that $\setS$ is isomorphic to the affine space $\F_3^3$.
For example, one can consider the theory of Hall triple systems \cite[IV.25]{Handbook}, or use a theorem of Kahn \cite{Kahn} (see also \cite[Theorem~7.4.5]{BB}) together with a counting argument.

\textbf{Characteristic \boldmath$\neq 2, 3$\unboldmath:}
Suppose $\setS$ contains at least two poor points $P_1$ and $P_2$.
Let $\Pi_i$ be a plane through $P_i$ such that $P_i$ has degree $3$ in $\Pi_i$.
Then $\Pi_1\neq\Pi_2$ (otherwise $\card{\Pi_i}\leq 7$, contradicting Lemma~\ref{smallcc}) and it is possible to choose a line $\ell_i$ in $\Pi_i$ through $P_i$, such that $\ell_1$ and $\ell_2$ are skew lines.
By Corollary~\ref{deg4}, for each $i=1,2$, $\card{\ell_i}=q_i+1$, where $q_i=p^{k_i}$ with $k_i\geq 1$, and $p=\kar{\D}\geq 5$.
By Proposition~\ref{threelines} all the points in $\ell_1\cup\ell_2\setminus\{P_1,P_2\}$ are bicoloured.
For any point $A\in\ell_2$, let $\Pi_A$ be the plane through $A$ and $\ell_1$.
Then by Lemma~\ref{lemma1} each line through $A$ in $\Pi_A$ has size $\geq 3$, except perhaps the line $P_1P_2$ (when $A=P_2$).
Since $\card{\ell_1}=q_1+1$, there are at least $q_1+1$ lines through $A$ in $\Pi_A$.
It follows that if $A\neq P_2$, then $\card{\Pi_A}\geq(3-1)\times(q_1+1)+1=2q_1+3$, with at least $q_1+2$ of these points not on $\ell_1$.
Similarly, if $A=P_2$, then $\card{\Pi_A}\geq 2q_1+2$ points, with at least $q_1+1$ not on $\ell_1$.
Therefore,
\[\card{\setS}\geq q_2(q_1+2)+q_1+1+\card{\ell_1}=q_1q_2+2q_1+2q_2+2.\]
If $p\geq 7$ or one of $k_i\geq2$, this gives $\card{\setS}> 51$.

Assume therefore that $q_1=q_2=p=5$.
Let $A\in\ell_2$ with $A\neq P_2$.
Considering $\Pi_A$ again, we have $\card{\Pi_A}\geq13$.
If $\card{\Pi_A}=13$, then all lines through $A$ have size $3$.
Since each point on $\ell_1$ except perhaps $P_1$ is bicoloured, it follows that all points in $\Pi_A$, except perhaps two on the line $AP_1$, are bicoloured, and it follows that $\Pi_A$ is an SGC.
However, no SG geometry on $13$ points has a line of size $6$ \cite{MR47:3207}.
Therefore, $\card{\Pi_A}\geq14$.
In the case $A=P_2$ we obtain $\card{\Pi_A}\geq12$.
Since now $\card{\ell_1}=6$, we obtain $\card{\setS}\geq5\times(14-6)+(12-6)+6>51$.

We have shown that if there is more than one poor point then $\card{\setS}> 51$.
We now assume that all points, except perhaps one, are rich.
First consider the case where there are no rich bicoloured points.
Then there is at most one bicoloured point (the poor point).
By removing this point, if necessary, we obtain an MR geometry $\setS'$.
Since all MR geometries up to $15$ points need characteristic $2$ or $3$ to embed \cite{PreSwa2}, we obtain that each plane $\Pi$ of $\setS'$ has at least $16$ points.
By putting the poor point back, we do not create any new planes, because it is bicoloured.
(A plane in $\setS$ not in $\setS'$ would be generated by a line of $\setS'$ and the bicoloured point, contradicting Lemma~\ref{lemma1}).
Therefore, all planes of $\setS$ have size $\geq 16$.
If all lines have size $\geq 5$, then by choosing any plane and any point outside the plane, and considering the lines through this point, we obtain $\card{\setS}\geq 4\times 16+1>51$.
Otherwise some line $\ell$ has size $\leq 4$.
Choose a rich point $P\in\ell$.
Choose any plane $\Pi$ through $P$ not containing $\ell$.
Because $P$ is rich, $P$ has degree $\geq 4$ in $\Pi$.
Then
\[\card{\setS}\geq 4\times(16-\card{\ell})+\card{\ell}\geq 4\times(16-4)+4>51.\]

We now assume that there exist rich bicoloured points.
If the contraction at some rich bicoloured point has size $\geq 25$, then by Lemma~\ref{lemma1} we have $\card{\setS}\geq25\times(3-1)+1=51$, and we are finished.
Assume therefore that the contraction at each rich bicoloured point has size $\leq 24$.
By Lemma~\ref{list} each such contraction must contain the geometry $20.2$, and be contained in $PG(2,5)$.
We now claim that each plane of $\setS$ has size $\geq 11$, and then conclude $\card{\setS}\geq51$.

Consider any plane $\Pi$ of $\setS$.
By Lemma~\ref{twolines} there is a rich point $P\notin\Pi$.
By Lemma~\ref{smallcc}, $\card{\Pi}\geq 9$.
If $\card{\Pi}=9$, it must be $AG(2,3)$, again by Lemma~\ref{smallcc}, giving that $AG(2,3)$ must be contained in the contraction $\setS/ P$.
However, $AG(2,3)$ is not embeddable in $PG(2,5)$ by Lemma~\ref{smallcc}.
This contradiction gives that any plane $\Pi$ has at least $10$ points.
However, the only $CC$ on $10$ points needs characteristic $3$ (again Lemma~\ref{smallcc}), and therefore, each plane has at least $11$ points.

Consider any rich, bicoloured point $P$.
Suppose that all the lines through $P$ corresponding to the points of $20.2$ in the contraction at $P$ have size $\geq 4$.
Then $\card{\setS}\geq20\times(4-1)+1>51$.
Therefore, assume that one of these lines $\ell$ has size $3$.
The degree of the point in $\setS/ P$ corresponding to $\ell$ is $6$, since the degree of each point of $20.2$ is $6$, and therefore, there are $6$ planes containing $\ell$.
This gives $\card{\setS}\geq(11-3)\times 6+3=51$, thereby finishing the proof.
\end{proof}

\section{Finite equivalents}\label{finite}
We now consider statements involving projective spaces over finite fields that are equivalent (in an elementary way) to Kelly's Theorem and the Hirzebruch Lemma.
These considerations lead to quantitative problems (Problem~\ref{problem3})
that would strengthen these two results.

Recall that in Section~\ref{proofthm3} we defined a $k$-SGC to be a subset $\setS$ of some projective space with at least $k$ points on each line determined by $\setS$.
Let $f_{k,n}(p)$ be the smallest size of an $n$-dimensional $k$-SGC in a projective space over $\F$, where $\F$ ranges over all fields of characteristic $p$.
If there does not exist such an $n$-dimensional $k$-SGC, set $f_{k,n}(p)=\infty$.
By Lemma~\ref{HN} it is sufficient to take this minimum over the finite fields $\F_{p^k}$.
To avoid trivial cases, we assume $k\geq 3$ and $n\geq 2$.
Thus the first non-trivial case is $f_{3,2}(p)$.
\begin{proposition}
\[ f_{3,2}(p)=\begin{cases} 7 & \text{if }p=2,\\ 9 & \text{for prime }p > 2. \end{cases}\]
\end{proposition}
\begin{proof}
We use the results of Kelly and Nwankpa \cite{MR47:3207} (see also Lemma~\ref{smallcc}).
Since a non-collinear SGC has at least $7$ points, we have $f_{3,2}(p)\geq 7$.
The only $7$-point non-collinear SGC is the Fano plane $PG(2,2)$, which embeds if and only if $\kar{\F}=2$.
An SGC with more than $7$ points must have at least $9$ points, with equality if and only if it is $AG(2,3)$.
This embeds if and only if $\kar{\F}=3$ or has an element of multiplicative order $3$.
However, if $p> 3$, then either $p\equiv 1\pmod{3}$, and then the multiplicative group of $\F_p$ has an element of order $3$, or $p\equiv 2\pmod{3}$, and then the multiplicative group of $\F_{p^2}$ has an element of order $3$.
This proves the proposition.
\end{proof}
The next two simplest cases are $f_{3,3}(p)$ and $f_{4,2}(p)$.
In the following proposition, by asserting that two statements are ``equivalent'', we mean that there is a known elementary proof that each implies the other.
\begin{proposition}
Kelly's Theorem is equivalent to the unboundedness of the set \[\{f_{3,3}(p): \text{$p$ is prime}\}.\]
The Hirzebruch Lemma is equivalent to the unboundedness of the set \[\{f_{4,2}(p): \text{$p$ is prime}\}.\]
\end{proposition}
\begin{proof}
We prove only the second statement. (The proof of the first statement is similar.)

First assume the Hirzebruch Lemma.
Suppose $f_{4,2}(p)\leq K$ for all primes $p$.
Since there are only finitely many two-dimensional $4$-SG geometries of size $\leq K$, there exists a single two-dimensional $4$-SG geometry $\setS$ embeddable in $P^2(\F)$ for fields $\F$ of infinitely many different characteristics.
The embeddability of $\setS$ into $P^2(\F)$ is a sentence $\fhi$ in the first order theory of fields.
Let
\[ \Sigma = \{\fhi\}\cup\{\underbrace{1+\dots+1}_{n}\neq 0 : n\geq 2\}\cup\{\text{field axioms}\}.\]
Any finite subset of $\Sigma$ is valid in some field $\F$ of sufficiently large characteristic $p$.
By the compactness theorem of first order logic there exists a model of $\Sigma$, i.e., a field $\F$ of characteristic $0$ into which $\setS$ embeds.
By Lemma~\ref{HN}, $\setS$ can be coordinatised by the algebraic closure of the prime field of $\F$, hence by $\C$.
Therefore, $\setS$ is embeddable in $P^2(\C)$, contradicting the Hirzebruch Lemma.

Conversely, assume that $f_{4,2}(p)$ is unbounded.
Suppose that the Hirzebruch Lemma is false.
Thus there exists a non-collinear $4$-SG geometry $\setS$ embeddable in $P^2(\C)$.
By Lemma~\ref{HN}, $\setS$ embeds in $P^2(\overline{\Q})$, where $\overline{\Q}$ is the algebraic closure of $\Q$.
Since $f_{4,2}(p)$ is unbounded, $\setS$ is not embeddable in $P^2(\F)$ for any $\F$ of characteristic $p$, for infinitely many primes $p$.
The non-embeddability of $\setS$ is again a first order sentence $\neg\fhi$.
Let $\psi_n$ be the first-order sentence asserting that any polynomial of degree $n$ has a root.
Let
\[ \Sigma=\{\neg\fhi\}\cup\{\underbrace{1+\dots+1}_{n}\neq 0 : n\geq 2\}\cup\{\psi_n: n\geq 2\}\cup\{\text{field axioms}\}.\]
Again it is clear that any finite subset of $\Sigma$ is valid in some field $\F$ of sufficiently large characteristic $p$.
The compactness theorem now gives an algebraically closed field $\F$ of characteristic $0$ such that $\setS$ does not embed in $P^2(\F)$.
However, since $\F$ contains $\overline{\Q}$, $\setS$ does embed, a contradiction.
\end{proof}

Since there is an elementary proof of Kelly's Theorem \cite{EPS}, we now have an elementary proof that $f_{3,3}(p)$ is unbounded.
Unfortunately this proof is not constructive, and gives no information on the size of $f_{3,3}(p)$.
See Problem~\ref{problem3} in Section~\ref{problems}.
 We know the following.
\begin{proposition}
\begin{equation} f_{3,3}(p)\begin{cases} =15, & p=2,\\ =27, & p=3, \\ \geq 51, & p>3,\\ \leq 3p^2, & p\geq 3.\end{cases}\label{eq1}
\end{equation}
\begin{equation} f_{4,2}(p)\begin{cases} =16, & p=2,\\ =13, & p=3,\\ =20, & p=5,\\ \geq 20, & p>3,\\ \leq 4p, & p>3.\end{cases}\label{eq2}
\end{equation}
\end{proposition}
\begin{proof}
The bound $f_{3,3}(p)\leq 3p^2$ follows by considering three parallel planes in the affine space $AG(3,p)$.
The remaining statements in \eqref{eq1} are implied by Theorem~\ref{thm3}.

The bound $f_{4,2}(p)\leq 4p$ follows by considering four parallel lines in the affine plane $AG(2,p)$, $p\geq 5$.
The remaining statements in \eqref{eq2} are implied by Corollary~\ref{4sgc}.
\end{proof}

Since the Hirzebruch Lemma is used in the proof of Kelly's Theorem, there should be a simple proof that the unboundedness of $f_{4,2}(p)$ implies the unboundedness of $f_{3,3}(p)$.
The following says a little more.
\begin{proposition}
$f_{3,3}(p)\geq\min\{2f_{4,2}(p)+1,6p+3\}$.
\end{proposition}
\begin{proof}
Let $\setS$ be a $3$-dimensional SGC in $P^3(\F_{p^k})$.
Let $P\in\setS$.
If the contraction of $\setS$ at $P$ is a $4$-SGC, then, since there are at least three points on each line through $P$, $\card{\setS}\geq 2f_{4,2}(p)+1$.
Otherwise there is a plane $\Pi$ through $P$ such that $\setS\cap\Pi$ is covered by three lines through $P$.
By Corollary~\ref{deg4} it follows that $\card{\setS\cap\Pi}\geq 3p+1$, and by considering the lines through a point of $\setS$ outside $\Pi$, we obtain $\card{\setS}\geq2\card{\setS\cap\Pi}+1\geq 6p+3$.
\end{proof}
Conversely, one could hope to prove the unboundedness of $f_{4,2}(p)$ directly from the unboundedness of $f_{3,3}(p)$, thus giving an elementary proof of the Hirzebruch Lemma.
See Problem~\ref{implication} in the next section.

\section{Open problems}\label{problems}
\subsection{Order and roots of unity}
As observed by Motzkin \cite{MR12:849c}, his proof of the SG Theorem  also holds for any ordered division ring, since the Pasch axiom still holds.
The same is true for the proof of Theorem~\ref{thm1} given in Section~\ref{introduction}.
Note that by the Artin-Schreier-Pickert-Szele Theorem (see \cite[Chapter~6]{Lam}), the ordered division rings are exactly the \emph{formally real} division rings (i.e.\ $-1$ is not a sum of products of squares).
\begin{corollary}
Any chromatic configuration in any projective $n$-space over a formally real division ring $\D$ is collinear.
\end{corollary}
In this regard, Kelly \cite{MR87k:14047} says:
\begin{center}
\emph{whether `Sylvester implies order' is an intriguing open question\dots}
\end{center}
\begin{problem}
Characterise the division rings \textup{(}or fields\textup{)} $\D$ over which
\begin{enumerate}
\item all SG configurations are collinear, or
\item all chromatic configurations are collinear.
\end{enumerate}
Are they exactly the formally real ones?
\end{problem}
Although there are examples of non-desarguesian projective planes that fail the axiom of Pasch even though all SGCs are collinear \cite{MR81c:05015, MR86a:51008}, we have no example of a division ring over which all SGCs are collinear that is not formally real as well.

Call a division ring \emph{root-free} if it does not contain an element of finite multiplicative order other than $\pm 1$.
It is easily observed that formally real division rings are root-free, and any root-free division ring either has characteristic $0$, or is a purely trancendental extension of $\F_2$.
Proposition~\ref{prop2} implies the following.
\begin{proposition}
If $\D$ is not root-free, then there exist non-collinear SGCs over $\D$.
\end{proposition}
By the results of Kelly and Nwankpa \cite{MR47:3207} it can be seen that if $\D$ is root-free and of characteristic $0$, then a $2$-dimensional SGC in $P^n(\D)$ must have at least $15$ points.
There exist root-free fields of characteristic $0$ admitting non-collinear SGCs, e.g., $\Q(\sqrt{-7})$ \cite{vW}.
This is the only known example of a non-collinear SGC over $\C$ that is not obtained from a root of unity via Proposition~\ref{prop2}.
The field in the following problem is also root-free.
\begin{problem}
Does there exist a non-collinear SGC over $\Q(\sqrt{-2})$?
\end{problem}

\subsection{Complex MR configurations}
Although  $P^2(\C)$ contains non-col\-lin\-ear SGCs (the smallest being the $9$-point configuration in the Introduction), we have no example of an MRC in $P^2(\C)$.
In fact, all known CCs in $P^n(\C)$ are SGCs.
\begin{problem}\label{p3a}
Does there exist a non-collinear MRC in $P^2(\C)$?
\end{problem}
It is known that if there exists a two-dimensional MRC in $P^n(\C)$ then it must have at least $19$ points \cite{PreSwa2, vW}.

\subsection{SG configurations over the quaternions}
The division ring of quaternions is the only non-commutative division ring for which an upper bound for the dimension of an SGC is known (the Quaternion Theorem stated in the Introduction).
\begin{problem}
Is an SGC in $P^n(\HQ)$ always coplanar?
\end{problem}
By Theorem~\ref{thm3}, a configuration giving a counterexample to the above question would have at least $51$ points.
It would also be interesting to find any non-commutative division ring of characteristic $0$ with a non-coplanar SGC. 

\subsection{SG and MR problems in characteristic \boldmath$2$\unboldmath}
Although the following discussion can be generalised to division rings, for simplicity we only consider fields in this section.
Let $\F$ be a field of finite characteristic $p$.
If $\F$ contains a finite field $\K$ of at least $3$ elements, then $P^2(\F)$ contains $P^2(\K)$, which is a non-collinear $4$-SGC and an MRC (since all finite projective planes except the Fano plane have blocking sets \cite[Example~8.1.3]{Batten}), and $\F^2$ contains $\K^2$, which is a non-collinear SGC.
Therefore, the only interesting fields of finite characteristic are those with only one finite subfield, of size $2$.
These are exactly the purely trancendental extensions of $\F_2$.
Since we work only with finite configurations, we may without loss assume that $\F$ has finite trancendence degree.
Furthermore, given any infinite field $\F$, any geometry embeddable in $\F(x)^2$ or $P^2(\F(x))$ can also be embedded in $\F^2$ or $P^2(\F)$, by substituting an appropriate element of $\F$ into $x$ in the coordinates of the configuration over the rational function field $\F(x)$.
It follows that we only have to consider $\F_2(x)$.
\begin{problem}\label{char2}\mbox{}
\begin{itemize}
\item Does there exist a non-collinear SGC \textup{(}or a non-collinear CC\textup{)} in $\F_2(x)^2$?
\item Does there exist a non-collinear MRC in $P^2(\F_2(x))$?
\item Does there exist a non-collinear $4$-SGC in $P^2(\F_2(x))$?
\end{itemize}
\end{problem}
It is known that a non-collinear SGC in $\F_2(x)^2$ has at least $15$ points (by the results in \cite{MR47:3207}), and a non-collinear MRC in $P^2(\F_2(x))$ needs at least $19$ points \cite{vW}.
By Lemma~\ref{list} a non-collinear $4$-SGC over $\F_2(x)$ needs at least $22$ points.

\subsection{Quantitative formulations of Kelly's Theorem and the Hirzebruch Lemma}
See Section~\ref{finite} for partial results on the following problem.
\begin{problem}\label{problem3}
Determine $f_{3,3}(p)$ and $f_{4,2}(p)$ in terms of $p$.
In particular, is $f_{3,3}(p)=3p^2$ for $p\geq 5$? Is $f_{4,2}(p)=4p$ for $p\geq 7$?
\end{problem}

\begin{problem}\label{implication}
Find a direct proof of the unboundedness of $f_{4,2}(p)$ from the unboundedness of $f_{3,3}(p)$.
\end{problem}

\end{document}